\documentclass{article}

\usepackage[colorlinks,bookmarksopen,bookmarksnumbered,citecolor=red,urlcolor=red]{hyperref}
\usepackage[numbered, framed, bw]{mcode} 
\usepackage{amsmath, amssymb, a4wide}

\usepackage{algorithmic}
\usepackage{algorithm}
\usepackage{tikz}
\usetikzlibrary{automata,positioning,shapes.geometric}
\usetikzlibrary{arrows.meta}
\tikzset{>={Latex[width=2mm,length=2mm]}}
\usepackage{circuitikz}

\newcommand{\vect}{\operatorname{vec}}

\algsetup{linenodelimiter=.}

\usepackage{setspace}
\let\Algorithm\algorithm
\renewcommand\algorithm[1][]{\Algorithm[#1]\setstretch{1.2}} 

\newtheorem{theorem}{Theorem}[section]

\newtheorem{problem}{Problem}[section]

\newcommand{\spann}{\operatorname{span}}

\begin{document}


\title{Numerical solution of Lyapunov equations related to\\ Markov jump linear systems}

\author{Tobias Damm \and Kazuhiro Sato \thanks{K. Sato is with Department of Systems Science, 
Graduate School of Informatics, 
Kyoto University,  Japan, email: kazuhirosato@i.kyoto-u.ac.jp} \and
  Axel Vierling
\thanks{T. Damm and A. Vierling are with the University of Kaiserslautern,
Department of Mathematics, 67663 Kaiserslautern, Germany, email:
damm@mathematik.uni-kl.de, axel4ling@googlemail.com}
}
\date{March 2017}
\maketitle
\begin{abstract}
  We suggest and compare different methods for the numerical solution of Lyapunov like equations with
application to control of Markovian jump linear systems. 
First, we consider fixed point iterations and associated
Krylov subspace formulations.
Second, we reformulate the equation as an optimization problem and
consider  steepest descent, conjugate gradient, and a trust-region method. 

Numerical experiments illustrate that for large-scale problems the
trust-region method is more effective than the steepest descent and
the conjugate gradient methods. The fixed-point approach, however,
is superior to the optimization methods. As an application we consider
a networked control system, where the Markov jumps are induced by the
wireless communication protocol.
\end{abstract}



\maketitle

\section{Introduction}
Markovian jump linear systems 
are expressed by
\begin{align}
\begin{cases}
\dot{x} = A(r(t))x +B(r(t))u,\\
y = C(r(t))x,
\end{cases} \label{Markov}
\end{align}
where $x\in\mathbb{R}^n$, $u\in\mathbb{R}^m$, and $y\in \mathbb{R}^p$
are the state, input, and output, respectively. 
The parameter $r(t)$ denotes a continuous-time Markov process on a
probability space, which takes values in a finite set
$S:=\{1,2,\ldots,N\}$ with transition probabilities given by 
\begin{align*}
{\rm Pr} \left( r(t+\delta)=j | r(t)=i\right) = 
\begin{cases}
\pi_{ij}\delta + o(\delta)\quad\quad\,\, {\rm if}\,\, i\neq j, \\
1+\pi_{ii}\delta + o(\delta)\,\,\,\, {\rm if}\,\, i=j,
\end{cases}
\end{align*}
where $\delta>0$, and $\pi_{ij}$ denotes the transition probability rate from mode $i$ to mode $j$ when $i\neq j$.
Furthermore, for all $i\in S$, $\pi_{ij}$ satisfies
$\pi_{ij} \geq 0$ $(i\neq j)$ and $\pi_{ii} = -\sum_{j\in S,\, j\neq i} \pi_{ij}$.
The Markov process $\{r(t),\, t>0\}$ is assumed to have an initial process $r(0)=(\mu_1,\mu_2,\ldots,\mu_N)$.
The matrices $A(r(t))$, $B(r(t))$, and $C(r(t))$ are contained in $\{ A_1,A_2,\ldots, A_N\}$, $\{ B_1,B_2,\ldots, B_N\}$, and $\{ C_1,C_2,\ldots, C_N\}$, respectively,
and if $r(t)=i\in S$, we have $A(r(t))=A_i$, $B(r(t))=B_i$, and $C(r(t))=C_i$.
Applications to electric power systems have been considered in
\cite{loparo1990probabilistic, ugrinovskii2005decentralized}, others
are mentioned e.g.\ in the survey paper \cite{shi2015survey} or the
monograph \cite{CostFrag13}. In subsection \ref{sec:WLAN} we sketch an
application to networked control system.

The system \eqref{Markov} is called mean-square stable if the solution $x(t)$ of
$ 
\dot{x}(t) = A(r(t))x(t)
$ 
satisfies
\begin{align*}
\lim_{t\rightarrow \infty} E(\|x(t)\|_2^2) =0
\end{align*}
for any initial condition $x(0)=x_0$ and initial distribution for $r(0)=r_0$, where $\|\cdot\|_2$ denotes the Euclidean norm and $E(\cdot)$ the expected value.
In \cite{costa1999continuous}, for mean-square stable Markovian jump linear systems \eqref{Markov}, the $H^2$ norm has been defined as
\begin{align*}
\|G\|_{H^2}:= \sqrt{\left( \sum_{s=1}^m \sum_{i=1}^N \mu_{i} \|y_{s,i}\|_2^2 \right) },
\end{align*}
where $\|y_{s,i}\|_2 := \sqrt{ \int_0^{\infty} E( y^T_{s,i}(t)y_{s,i}(t)) dt}$ and $y_{s,i}$ is the output $\{y(t)|t>0\}$ when
\begin{itemize}
\item the input is given by $u(t)=e_s \delta(t)$, where $\delta (t)$ is the unit impulse, and $e_s$ is the $m$ dimensional unit vector formed by having $1$ at the $s$th position and zero elsewhere.
\item $x(0)=0$ and $r(0)=i$.
\end{itemize}
For $N=1$, the definition reduces to the usual $H^2$ norm.
If
\eqref{Markov} is mean-square stable, then by \cite{costa1999continuous}
\begin{align*}
\|G\|_{H^2}^2 =\sum_{i=1}^N {\rm tr}(C_iP_iC_i^T)= \sum_{i=1}^N \mu_i {\rm tr}(B_i^T Q_i B_i),
\end{align*}
where $(P_1,P_2,\ldots,P_N)$ and $(Q_1,Q_2,\ldots,Q_N)$ are the unique
solutions of the coupled Lyapunov equations
\begin{align}
& A_i P_i+P_iA_i^T+\sum_{j=1}^N \pi_{ji} P_j+\mu_iB_iB_i^T = 0, \label{lya1} \\
& A_i^T Q_i+Q_iA_i+\sum_{j=1}^N \pi_{ij} Q_j+C_i^TC_i = 0 \label{lya2}
\end{align}
for $i=1,2,\ldots,N$, respectively.
Using the above $H^2$ norm concept, \cite{costa1999continuous} and \cite{sun2012optimal} have studied an $H^2$ optimal state-feedback control and an $H^2$ optimal model reduction method for Markovian jump linear systems \eqref{Markov}, respectively.
Thus it is important to study effective algorithms for solving the Lyapunov like equations \eqref{lya1} and \eqref{lya2}.
In a unified form we look for solutions $X=(X_1,\ldots,X_N)$ of
coupled equations 
\begin{align}
0&=f_i(X)  := A_iX_i +X_i A_i^T + \sum_{j=1}^N \gamma_{ij} X_j + Y_i\quad (i=1,2,\ldots,N)\;. \label{1}
\end{align}
Here $A_i, Y_i=Y_i^T \in \mathbb{R}^{n\times n}$ are constant
matrices, and $\gamma_{ij} \geq 0$ $(i\neq j)$ and $\gamma_{ii} =
-\sum_{j\in S,\, j\neq i} \gamma_{ij}$, or more generally, just $\gamma_{ii}<0$.
Eqs.\,\eqref{lya1} and \eqref{lya2} are special cases of \eqref{1}.
If $\gamma_{ij}=0$ for all $i\neq j$, then \eqref{1}  just consists of uncoupled standard Lyapunov equations.
For the standard Lyapunov equation, several algorithms have been developed,
see \cite{damm2008direct, kjelgaard2009numerical}, and references
therein. Coupled Lyapunov equations of type \eqref{1}  have been
solved e.g.\ in
\cite{Born95, CollHode97, DingChen06, WangLam08, ZhouLam08, LiZhou11},
partly by fixed-point iterations, partly by optimization methods.
However, all methods were only applied to examples of small
dimensions. In trying to extend methods of model order reduction to Markov
jump linear systems, we found the necessity to develop a more
efficient solver. In this paper we compare different classes of algorithms for their
applicabilty at least to medium-sized problems with $nN\approx 1000$ or larger.

Our first and main approach follows the ideas in \cite{damm2008direct} and
uses a fixed point iteration. Here the efficiency can be improved by
considering Krylov-subspace methods and appropriate vectorization. The
idea is rather simple, but works much better than other more intricate
methods. 

In our second approach, we reformulate the equation as an optimization problem on the product
space of $N$ Euclidean matrix spaces and derive the gradient and the Hessian of the objective function.
The gradient is used to develop a steepest descent and a conjugate
gradient method; the Hessian is applied to establish a trust-region method.

Numerical experiments illustrate that all our methods work, but 
only the fixed-point iteration lends itself for larger problems.

\section{Preliminaries}



Let $\mathcal{H}^n$ denote the space of symmetric (i.e.\ real Hermitian) $n\times n$
matrices and $\mathcal{H}_+^n$ the cone of nonnegative definite
matrices.

Then we set
\begin{align*}
\mathcal{H}=\underbrace{\mathcal{H}^n\times\cdots\times\mathcal{H}^n}_{N\text{-times}}\quad
\text{ and }\quad\mathcal{H}_+=\mathcal{H}_+^n\times\cdots\times\mathcal{H}_+^n\;,
\end{align*}
such that $\mathcal{H}$ is an ordered real vector space with ordering
cone $\mathcal{H}_+$. For $A\in\mathbb{R}^{n\times n}$, we define the
Lyapunov operator
\begin{align*}
\mathcal{L}_A:\mathcal{H}^n\to\mathcal{H}^n\quad\text{  by }\quad
\mathcal{L}_A(H)=AH+HA^T\;.
\end{align*}
On $\mathcal{H}$ we consider the blockwise Lyapunov operator
$\mathcal{L}:\mathcal{H}\to\mathcal{H}$, defined by
\begin{align*}
\mathcal{L}(X_1,\ldots,X_N)=(Z_1,\ldots,Z_N)\quad\text{ with }\quad
  Z_i=\mathcal{L}_{A_i}(X_i)+\gamma_{ii}
  X_i=\mathcal{L}_{A_i+\tfrac{\gamma_{ii}}2I}(X_i) 
\end{align*}
and the positive operator
$\Pi:\mathcal{H}\to\mathcal{H}$, $\Pi(\mathcal{H}_+)\subset\mathcal{H}_+$, defined by
\begin{align*}
\Pi(X_1,\ldots,X_N)=(Z_1,\ldots,Z_N)\quad\text{ with }\quad Z_i=\sum_{i\neq
  j}\gamma_{ij}X_j\;.
\end{align*}
From \cite[Thm. 3.15]{CostFrag13} and
also \cite{Schn65, Damm04} we cite a well-known stability result.
\begin{theorem}
The following statements are equivalent.
  \begin{itemize}
  \item[(a)] System \eqref{Markov} is asymptotically mean-square
    stable.
  \item[(b)]
    $\sigma(\mathcal{L}+\Pi)\subset\mathbb{C}_-=\left\{\lambda\in\mathbb{C}\;\big|\;\Re
      \lambda<0\right\}$.
  \item[(c)] $\sigma(\mathcal{L})\subset\mathbb{C}_-$ and
$\rho(\mathcal{L}^{-1}\Pi)<1$.
\item[(d)] $\exists X\in\mathcal{H}_+$: $(\mathcal{L}+\Pi)(X)<0$.
\item[(e)] $\forall Y \in\mathcal{H}_+$:$\exists X\in\mathcal{H}_+$: $(\mathcal{L}+\Pi)(X)=-Y$.
  \end{itemize}
\end{theorem}
\noindent
The coupled Lyapunov equations
\eqref{1}  can be written in the form
\begin{align}\label{eq:lyapXY}
 ( \mathcal{L}+\Pi) (X)=-Y\;,
\end{align}
where $Y\in\mathcal{H}$. Under the assumptions of asymptotic
mean-square stability there exists a
unique solution $X\in\mathcal{H}$. Moreover, if $Y\in\mathcal{H}_+$,
then $X\in\mathcal{H}_+$. Since $X$ contains
$\frac12{Nn(n+1)}$ scalar unknowns, a direct solution e.g.\ via
Kronecker-product representation and Gaussian elimination has
complexity $O(N^3n^6)$.

However, it is well-known, that a single Lyapunov equation for an
unknown $n\times n$-matrix can be solved with $O(n^3)$ operations by
the Bartels-Stewart-algorithm \cite{BartStew72}. Since
$\mathcal{L}^{-1}(X)$ for  $X\in\mathcal{H}$ is obtained by solving
$N$ independent Lyapunov equations, the cost of evaluating
$\mathcal{L}^{-1}(X)$ is only $O(Nn^3)$.  Based on this observation we
suggest some fixed point iterations.


\section{Fixed point formulations and Krylov subspace methods}



\subsection{Jacobi and Gauss-Seidel schemes}

\label{sec:fixed-point-form}
If $\rho(\mathcal{L}^{-1}\Pi)<1$ then the solution to \eqref{eq:lyapXY}  is obtained as the limit of
the iterative scheme
\begin{align*}
  X^{(k+1)}&=-\mathcal{L}^{-1}\left(\Pi(X^{(k)})+Y\right)\;.
\end{align*}
Blockwise, we get the Jacobi-type fixed point iteration
\begin{align*}
  X_1^{(k+1)}&=-\mathcal{L}_{A_1+\tfrac{\gamma_{11}}2I}^{-1}\left(Y_1+\gamma_{12}X_2^{(k)}+\gamma_{13}X_3^{(k)}+\ldots+\gamma_{1N}X_N^{(k)}\right)\\
&\;\;\vdots\\
X_N^{(k+1)}&=-\mathcal{L}_{A_N+\tfrac{\gamma_{NN}}2I}^{-1}\left(Y_N+\gamma_{N1}X_1^{(k)}+\gamma_{N2}X_2^{(k)}+\ldots+\gamma_{N,N-1}X_{N-1}^{(k)}\right)\;.
\end{align*}
In the $i$-th row, it is natural to replace $X_j^{(k)}$ by
$X_j^{(k+1)}$ for $j<i$, because this update is already available.
Thus we get the Gauss-Seidel-type scheme
\begin{align*}
  X_1^{(k+1)}&=-\mathcal{L}_{A_1+\tfrac{\gamma_{11}}2I}^{-1}\left(Y_1+\gamma_{12}X_2^{(k)}+\gamma_{13}X_3^{(k)}+\ldots+\gamma_{1N}X_N^{(k)}\right)\\
X_2^{(k+1)}&=-\mathcal{L}_{A_2+\tfrac{\gamma_{22}}2I}^{-1}\left(Y_2+\gamma_{21}X_1^{(k+1)}+\gamma_{23}X_3^{(k)}+\ldots+\gamma_{2N}X_N^{(k)}\right)\\
&\;\;\vdots\\
X_N^{(k+1)}&=-\mathcal{L}_{A_N+\tfrac{\gamma_{NN}}2I}^{-1}\left(Y_N+\gamma_{N1}X_1^{(k+1)}+\gamma_{N2}X_2^{(k+1)}+\ldots+\gamma_{N,N-1}X_{N-1}^{(k+1)}\right)\;.
\end{align*}

\subsection{Preconditioned Krylov subspace iterations}
\label{sec:prec-kryl-subsp}

Both schemes can be written in the form $X^{(k+1)}=T(X^{(k)})$ and interpreted as preconditioners for Krylov subspace
iterations as has been explained e.g.\ in \cite[Sec
4.2]{damm2008direct}. The basic idea is to find an optimal
approximation to the solution within the Krylov subspace $X^{(0)}+\spann\{X^{(1)},\ldots,X^{(k)}\}$.
 In the Jacobi formulation, this means that we set
$T_J=-\mathcal{L}^{-1}\Pi$ and apply
some standard Krylov subspace method to the equation
\begin{align*}
  (I-T_J)(X)=-\mathcal{L}^{-1}(Y)=\tilde Y\;.
\end{align*}
That is, we replace $-Y$ by $\tilde Y=-\mathcal{L}^{-1}(Y)$ and then solve with
the linear mapping $I-\mathcal{L}^{-1}\Pi$. 
In the Gauss-Seidel formulation, the linear mapping
$T_{\text{GS}}:X\to\tilde X$ is described by the scheme
\begin{align}\label{eq:defTGS}
  \begin{split}
    \tilde X_1&=-\mathcal{L}_{A_1+\tfrac{\gamma_{11}}2I}^{-1}\left(\gamma_{12}X_2+\gamma_{13}X_3+\ldots+\gamma_{1N}X_N\right)\\
    \tilde
    X_2&=-\mathcal{L}_{A_2+\tfrac{\gamma_{22}}2I}^{-1}\left(\gamma_{21}\tilde
      X_1+\gamma_{23}X_3+\ldots+\gamma_{2N}X_N\right)\\
    &\;\;\vdots\\
    \tilde
    X_N&=-\mathcal{L}_{A_N+\tfrac{\gamma_{NN}}2I}^{-1}\left(\gamma_{N1}\tilde
      X_1+\gamma_{N2}\tilde X_2+\ldots+\gamma_{N,N-1}\tilde
      X_{N-1}\right)\;.
  \end{split}
\end{align}
The update of the right hand side $Y\mapsto \tilde Y$ is obtained via
\begin{align}\label{eq:defYGS}
  \begin{split}
    \tilde Y_1&=-\mathcal{L}_{A_{1}-\tfrac{\gamma_{11}}2I}^{-1}(Y_1)\\
    \tilde Y_2 &=-\mathcal{L}_{A_{2}-\tfrac{\gamma_{22}}2I}^{-1}
    (Y_2+\gamma_{21}\tilde
    Y_1)\\
    &\;\;\vdots\\
    \tilde Y_N&=-\mathcal{L}_{A_{N}-\tfrac{\gamma_{NN}}2I}^{-1}
    (Y_N+\gamma_{N1}\tilde Y_1+\ldots+\gamma_{N,N-1}\tilde Y_{N-1})\;.
  \end{split}
\end{align}
Thus we can apply a Krylov subspace method to the equation
$(I-T_{\text{GS}})(X)=\tilde Y$.


\begin{algorithm}                      
\caption{Krylov-subspace method with Gauss-Seidel-type preconditioning
  for equation \eqref{eq:lyapXY}}         
\label{algorithm1}                          
\begin{algorithmic}[1]
\STATE Choose tolerance level \texttt{tol}
\STATE Define function $X\mapsto T_{GS}(X)$ according to
\eqref{eq:defTGS}
\STATE Compute preconditioned right-hand side $\tilde Y$ according to
\eqref{eq:defYGS}
\STATE Compute solution $X$ by Krylov subspace method applied to
$(I-T_{GS})(X)=\tilde Y$.
\end{algorithmic}
\end{algorithm}

\subsection{Avoiding loops by vectorization}
\label{sec:avoid-loops-vect}

Terms of the form
$\sum_{j=1}^N \gamma_{ij}X_j$  
can be vectorized efficiently. If $\vect X_i$ as before denotes the
vector obtained by stacking all columns of $X_i$ one above the other,
and (by a slight abuse of notation) $\vect X=[\vect X_1,\ldots,\vect X_N]$, then 
\begin{align} \label{eq:sum_vect}
  \vect\Big(\sum_{j=1}^N \gamma_{ij}X_j\Big)&=\vect X \Gamma_i^T\;,
\end{align}
where $\Gamma_i$ is the $i$-th row of $\Gamma$. Thus, we can write
the sum in \eqref{eq:sum_vect} as a matrix-vector product, which is
processed faster than a loop over $j$.

 \subsection{An implementation}
 \label{sec:implementation}

Some of the specific vectorization ideas can be seen more clearly in the
MATLAB\textsuperscript{\textregistered}-listing which we add for convenience. 
\begin{lstlisting}
function X=MJLSlyap(A,Y,G,tol)
% Solve algebraic Lyapunov equation for Markov jump linear system
% The n^2 x N matrices A and Y contain A1(:),...,AN(:) and Y1(:),...,YN(:)
% as columns. The matrix G=Gamma is N x N
% Krylov subspace approach with Gauss-Seidel-type preconditioning

N=size(G,2);n2=size(A,1);n=sqrt(n2);
d=diag(G);G=G'-diag(d);d=d/2;
A=reshape(A+reshape(eye(n),n2,1)*d',n,n,N);
Y=reshape(Y,n,n,N);
for i=1:N
    Yi=Y(:,:,i)+reshape(reshape(Y(:,:,1:i),n2,i)*G(1:i,i),n,n);
    Y(:,:,i)=lyap(A(:,:,i),Yi);
end  
X=reshape(bicgstab(@T,Y(:),tol),n,n,N);
    function Z = T(X)
        X=reshape(X,n,n,N);Z=zeros(size(X));
        for j=1:N
            Xj=lyap(A(:,:,j),reshape(reshape(X,n2,N)*G(:,j),n,n));
            Z(:,:,j)=X(:,:,j)-Xj;
            X(:,:,j)=Xj;
        end
        Z=Z(:);
    end
end
\end{lstlisting}

\section{Optimization-based approach} \label{sec2}

This section formulates an optimization problem for solving the
Lyapunov like equations \eqref{1} and suggests a number of methods to
solve this problem.

\subsection{Reformulation as optimization problem}

To solve \eqref{1},
we consider the following optimization problem.

\begin{problem}\label{problem1}
  \begin{align*}
    &{\rm minimize} \quad f(X) := \sum_{i=1}^N \|f_i(X_1,X_2,\ldots,X_N)\|_F^2 \\
    &{\rm subject\, to} \quad\displaystyle X\in M:= \mathbb{R}^{n\times n} \times \mathbb{R}^{n\times n}\times \cdots \times \mathbb{R}^{n\times n}.
  \end{align*}
\end{problem}


If $f(X)\approx 0$, we obtain $f_i(X)\approx 0$ $(i=1,2,\ldots,N)$, i.e., $X$ is an approximate solution of \eqref{1}.



\subsection{Gradient and Hessian of the objective function $f$} \label{sec3}

To develop optimization algorithms for Problem \ref{problem1}, this section derives the gradient of the objective function $f$, and then gives the Hessian of $f$.

The Fr\'echet derivative of $f$ at $X=(X_1,X_2,\ldots,X_N)$ in the direction $X'=(X'_1,X'_2,\ldots,X'_N)$ can be calculated as
\begin{align}
 {\rm D}f(X)[(X')]  
&=  2\sum_{i=1}^N {\rm tr} \left( (A_iX'_i+X'_iA_i^T +\sum_{j=1}^N \gamma_{ij} X'_j)f_i \right) = \sum_{i=1}^N{\rm tr} \left( X'^T_i D_i\right)  \label{8}
\end{align}
for
$
\displaystyle D_i:= 2 \Big(A_i^Tf_i+f_iA_i+\sum_{j=1}^N\gamma_{ji} f_j \Big).
$

\noindent
Since the gradient ${\rm grad} f(X)$ satisfies
$\displaystyle
 {\rm D}f(X)[(X')]  
=  \sum_{i=1}^N {\rm tr} (X'^T_i ({\rm grad} f(X))_{X_i})$, equation 
\eqref{8} yields
\begin{align}
 {\rm grad}\, f(X)  = ( D_1, D_2, \ldots, D_N). \label{Euclid_grad}
\end{align}

\noindent
Furthermore, from \eqref{1} and \eqref{Euclid_grad}, the Hessian of $f$ is given by
\begin{align}
{\rm Hess}\, f(X)[(\xi_1,\xi_2,\ldots,\xi_N)] 
 :={\rm D} {\rm grad}\, f (X)[ (\xi_1,\xi_2,\ldots,\xi_N)] 
= (D'_1,D'_2,\ldots,D'_N ), \label{Hess_f}
\end{align}
where
$\displaystyle D'_i = 2\Big( A_i^Tf'_i+f'_iA_i + \sum_{j=1}^N\gamma_{ji} f'_j \Big)$
with
$
f'_i= A_i\xi_i+\xi_iA_i^T + \sum_{j=1}^N \gamma_{ij} \xi_j$.

\subsection{Optimization algorithms for Problem \ref{problem1}} \label{sec4}

In the vector space $M$, optimization methods based on line search can be developed.
In such methods, given the current point $p_k\in M$, the search direction $d_k\in T_{p_k} M \cong M$, and the step size $t_k>0$, the next point $p_{k+1}\in M$ is computed as
\begin{align*}
p_{k+1} = p_k + t_k d_k.
\end{align*}
We consider three optimization algorithms for Problem \ref{problem1}.
To this end, 
for any $(\xi_1,\xi_2,\ldots, \xi_N), (\eta_1,\eta_2,\ldots,\eta_N) \in T_p M \cong M$, we define the inner product as
\begin{align*}
\langle (\xi_1,\xi_2,\ldots, \xi_N), (\eta_1,\eta_2,\ldots,\eta_N) \rangle:= \sum_{i=1}^N{\rm tr}(\xi_i^T\eta_i),
\end{align*}
and the induced norm by
$
\|(\xi_1,\xi_2,\ldots, \xi_N)\|:= \sqrt{ \langle (\xi_1,\xi_2,\ldots, \xi_N), (\xi_1,\xi_2,\ldots,\xi_N) \rangle}.
$

\subsubsection{Steepest descent method for Problem \ref{problem1}}


In the steepest descent method,
the negative gradient of the objective function $f$ at a current iterate $p_k\in M$
can be chosen as a search direction $d_k \in T_{x_k} M$ at $p_k$, i.e., $d_k:= -{\rm grad}\, f(p_k)$.
As a step size $t_k$, the following Armijo step size is popular \cite{nocedal2006numerical}:
Given a  point $p\in M$, a tangent vector $d\in T_p M$, and scalars $\bar{\alpha}>0, \beta, \sigma \in (0,1)$, 
the Armijo step size $t^A :=\beta^{\gamma} \bar{\alpha}$ is defined in such a way that $\gamma$ is the smallest nonnegative integer satisfying
\begin{align}
 f(p + \beta^{\gamma}\bar{\alpha} d) \leq f(p)+\sigma \langle {\rm grad}\, f(p), \beta^{\gamma}\bar{\alpha} d \rangle \label{Armijo}
\end{align}

\noindent
Algorithm \ref{algorithm1} describes the steepest descent method for solving Problem \ref{problem1}.

\begin{algorithm}                      
\caption{Steepest descent method for Problem \ref{problem1}.}         
\label{algorithm1}                          
\begin{algorithmic}[1]
\STATE Choose an initial point $p_0 \in M$.
\FOR{$k=0,1,2,\ldots$ }
\STATE Compute the search direction $d_k\in T_{p_k} M$ by
$d_k = -{\rm grad}\, f(p_k)$.
\STATE Compute the Armijo step size $t^A_k>0$ satisfying \eqref{Armijo}.
\STATE Compute the next iterate 
$p_{k+1} = p_k + t^A_k d_k$.
\ENDFOR
\end{algorithmic}
\end{algorithm}


\subsubsection{Conjugate gradient method for Problem \ref{problem1}}


In the conjugate gradient method, the search direction $d_k$ at the current point $p_k$ is computed as
\begin{align}
d_k = -{\rm grad}\, f(p_k) + \beta_k d_{k-1}, \label{conjugate}
\end{align}
where $\beta_k>0$.
The Dai and Yuan type parameter $\beta_k$ is given by 
\begin{align}
\beta_k = \frac{ \| {\rm grad}\, f(p_k)\|^2 }{ \langle d_{k-1}, {\rm grad}\, f(p_k)-{\rm grad}\, f(p_{k-1}) \rangle}. \label{dai}
\end{align}
In more detail, see \cite{dai1999nonlinear, sato2016dai}.
To guarantee the convergence property,
we use the step size $t_k$ satisfying 
\begin{align}
\begin{cases}
f(p_k+t_kd_k) \leq f(p_k) + c_1 t_k \langle {\rm grad}\, f(p_k), d_k \rangle, \\
\langle {\rm grad}\,f(p_k+t_k d_k), d_k \rangle \geq c_2 \langle {\rm grad}\, f(p_k), d_k \rangle,
\end{cases} \label{wolfe}
\end{align}
where $0<c_1<c_2<1$.
The condition \eqref{wolfe} is called the Wolfe condition.
Algorithm \ref{algorithm2} describes the conjugate gradient method for solving Problem \ref{problem1}.

\begin{algorithm}                      
\caption{Conjugate gradient method for Problem \ref{problem1}.}         
\label{algorithm2}                          
\begin{algorithmic}[1]
\STATE Choose an initial point $p_0 \in M$.
\STATE Set $d_0 = -{\rm grad}\, f(p_0)$.
\FOR{$k=0,1,2,\ldots$ }
\STATE Compute the step size $t^W_k>0$ satisfying \eqref{wolfe}.
\STATE Compute the next iterate $
p_{k+1} = p_k + t^W_k d_k.$
\STATE Set $\beta_{k+1}$ by \eqref{dai}.
\STATE Set $d_{k+1}$ by \eqref{conjugate}.
\ENDFOR
\end{algorithmic}
\end{algorithm}

\subsubsection{Trust-region method for Problem \ref{problem1}}


At each iterate $p$ in the trust-region method on the vector space $M$, 
we evaluate the quadratic model $\hat{m}_{p}$ of the objective function $f$ within a trust region:
\begin{align*}
 \hat{m}_{p}(d)  = f(p) + \langle {\rm grad}\,f(p), d \rangle +\frac{1}{2}  \langle {\rm Hess}\, f(p)[d], d \rangle
\end{align*}
for $d\in T_p M$,
where ${\rm grad}\,f(p)$ and ${\rm Hess}\, f(p)[d]$ are given by \eqref{Euclid_grad} and \eqref{Hess_f}, respectively.
A trust-region with a radius $\Delta>0$ at $p\in M$ is defined as a ball with center $0$ in $T_p M$.
Thus the trust-region subproblem at $p\in M$ with a radius $\Delta$ is defined as a problem of minimizing $\hat{m}_{p}(d)$ subject to $d\in T_p M$, $\|d\| \leq \Delta$.
This subproblem can be solved by the truncated conjugate gradient method \cite{absil2007trust}.
Then we compute the ratio of the decreases in the objective function $f$ and the model $\hat{m}_p$ attained by the resulting $d_*$ to decide
whether $d_*$ should be accepted and whether the trust-region with the radius $\Delta$ is appropriate.
Algorithm \ref{algorithm} describes the process.
The constants $\frac{1}{4}$ and $\frac{3}{4}$ in the condition expressions in Algorithm \ref{algorithm} are commonly used in the trust-region method for a general unconstrained optimization problem.
These values ensure the convergence properties of the algorithm \cite{absil2007trust}.

\begin{algorithm}                      
\caption{Trust-region method for Problem \ref{problem1}.}    \label{algorithm}     
\label{alg1}                          
\begin{algorithmic}[1]
\STATE Choose an initial point $p_0 \in M$ and parameters $\bar{\Delta}>0$, $\Delta_0\in (0,\bar{\Delta})$, $\rho'\in [0,\frac{1}{4})$.
\FOR{$k=0,1,2,\ldots$ }
\STATE Solve the following trust-region subproblem for $d$ to obtain $d_k\in T_{p_k} M$:\\
$
{\rm minimize}\quad \hat{m}_{p_k}(d) \quad {\rm subject\, to}\quad
  \|d\| \leq \Delta_k, \quad {\rm where}\quad  d \in T_{p_k}M \cong M.
$
\STATE Evaluate
$
\displaystyle\rho_k := \frac{ f(p_k) -f(p_k+d_k)}{ \hat{m}_{p_k}(0)- \hat{m}_{p_k} (d_k)}. 
$
\IF {$\rho_k<\frac{1}{4}$}
\STATE 
$\Delta_{k+1}=\frac{1}{4}\Delta_k$.
\ELSIF {$\rho_k>\frac{3}{4}$ and $\|d_k\| = \Delta_k$}
\STATE
$\Delta_{k+1} = \min (2\Delta_k,\bar{\Delta})$.
\ELSE 
\STATE
$\Delta_{k+1} = \Delta_k$.
\ENDIF
\IF {$\rho_k>\rho'$}
\STATE
$p_{k+1} = p_k + d_k$.
\ELSE
\STATE
$p_{k+1} = p_k$.
\ENDIF
\ENDFOR
\end{algorithmic}
\end{algorithm}

\subsection{Stopping criterion for Algorithms \ref{algorithm1}, \ref{algorithm2}, and \ref{algorithm}}

In practice, we need a stopping criterion for Algorithms \ref{algorithm1}, \ref{algorithm2}, and \ref{algorithm}.
In this paper, we stop the algorithm when $\|{\rm grad}\,f(p_k)\| < 10^{-5}$ or the iteration number reached 30000.

\section{Numerical experiments and an application} \label{sec5}
We first discuss some artificial and random examples to compare the
performance of Algorithms 1--4. It turns out that Algorithm 1 is our
method of choice.
Then we provide some details on an
application of Markov jump linear systems to networked control
systems. Our method works fine for this example as well.
All computations were carried out using MATLAB\textsuperscript{\textregistered} R2016b on an 
Intel\textsuperscript{\textregistered} Core(TM) i7-7500U CPU @  2.70GHz 2.90GHz and 16.0 GB RAM.
We report on  computation times, numbers of required iterations,
residuals and -- where available -- the absolute errors. Note that for
Algorithm 1 each evaluation of $T_{\text{GS}}$ is half an iteration
step. The residual is given as the square root of the objective
function $f$.

\subsection{A system with known solution}
Our first example is constructed with known solution. Let $\begin{pmatrix}
    \gamma_{11}&\gamma_{12}\\\gamma_{21}&\gamma_{22}
  \end{pmatrix}=\begin{pmatrix}
    -1&1\\2&-2
 \end{pmatrix}$ and
\begin{align*} 
A_1 &:= \begin{pmatrix}
-6&  4& -7&  6\\
  8& -4&-10& 10\\
 14&  6&  1&  7\\
   -21&-10& -6&-13\\
\end{pmatrix}, \quad
A_2 := \begin{pmatrix}
   -16&  4&  7& -1\\
     5&-17& -8& -2\\
    -2&  3&-19& -4\\
     4& 10& 25& -9
\end{pmatrix}, \\
 X_1 &:= \begin{pmatrix}
1&  1&  1&  1\\
1&  1&  1&  1\\
1&  1&  1&  1\\
1&  1&  1&  1
\end{pmatrix}, \quad
X_2 := \begin{pmatrix}
2&  1& -1& -2\\
1&  1&  0& -1\\
-1&  0&  1&  1\\
-2& -1&  1&  2
\end{pmatrix}.
\end{align*}
Then $Y_j=-(\mathcal{L}+\Pi)(X)$ is given by
\begin{align*}
Y_1 &:= \begin{pmatrix}
 5& - 1& -23&56\\
 -1&  -8& -31&48\\
-23& -31& -56&22\\
56&48&22&99
\end{pmatrix}, \quad
Y_2 := \begin{pmatrix}
 68& 6& -52& -50\\
 6&20& 8& -22\\
 -52& 8&42&14\\
 -50& -22&14&24
\end{pmatrix}\;. 
\end{align*}

For Algorithm 1, 2, 3, and 4, Table \ref{table1} shows computational times, iteration numbers, the
residual and the error
$\|\Delta\|_F$, i.e. the deviation of the
computed solution from the known solution, in the Frobenius norm.
According to Table \ref{table1}, Algorithm 1 is superior to the others.

\begin{table}[h]
\caption{Constructed example with $n=4$, $N=2$ and known solution} \label{table1}
  \begin{center}
    \begin{tabular}{|l|c|c|c|c|} \hline
            & Time (seconds) & Iteration & Residual& Error\\ \hline 
      Algorithm 1 & 0.0041 & 3 & 1.1e$-$11
&4.0e$-$13 \\
      Algorithm 2 & 2.5 & 1116 & 2.9e$-$05 
&7.5e$-$06 \\
	Algorithm 3 & 0.93 & 359 & 1.9e$-$07
&2.5e$-$08 \\
      Algorithm 4 & 0.19 & 14 &5.5e$-$08
&2.6e$-$09 \\ \hline
    \end{tabular}
  \end{center}
\end{table}

\subsection{Random systems}
We consider random  mean square
asymptotically systems.  In  our first experiments we always set $\Gamma=\begin{pmatrix}
    -0.3 & 0.3 \\
0.5&-0.5
 \end{pmatrix}$ and  generate $A_1, A_2, Y_1, Y_2\in
 \mathbb{R}^{n\times n}$ such that $-A_1,-A_2,Y_1,Y_2$ are symmetric
 positive definite. Under this assumption we have
 $(\mathcal{L}+\Pi)([I,I])=-2([A_1,A_2])<0$ implying asymptotic
 mean-square stability.

Table \ref{table2} shows the computational times, iteration numbers,
and the residuals $\sqrt{f(X_1,X_2)}$ in Algorithms 1,
2,  3, and 4 for different $n$. For $n=100$ Algorithms 2 and 3
turn out to be inadequate and for larger dimensions Algorithm 4
becomes too slow.
\begin{table}[h]
\caption{Random examples with fixed $\Gamma\in\mathbb{R}^{2\times 2}$
  and positive definite $A_1, A_2$} \label{table2}
  \begin{center}
    \begin{tabular}{|l|c|c|c|c|} \hline
 & $n$& Time (seconds) & Iteration & Residual \\ \hline 
      Algorithm 1 & 100&0.11 & 5 & 7.1e$-$10 \\
	Algorithm 2 & 100&113 & 30000 & 1.0e02\\
      Algorithm 3 & 100&131 & 30000 & 5.9e00 \\
	Algorithm 4 & 100&33 & 30 & 6.9e$-$08 \\\hline
      Algorithm 1 &300& 1.3 & 4 & 1.1e$-$09 \\
	Algorithm 4 &300& 745 & 35 & 6.4e$-$08\\
      \hline
      Algorithm 1 & 1200&88.0 & 3.5 & 5.5e$-$08\\\hline
    \end{tabular}
  \end{center}
\end{table}


In the next experiments we generate random nonsymmetric matrices
$A_j\in\mathbb{R}^{n\times n}$, $j=1,\ldots,N$ with
$\sigma(A_j)\subset\mathbb{C}_-$ and 
$\Gamma=(\gamma_{ij})\in\mathbb{R}^{N\times N}$ with $\gamma_{ii}<0$
and $\gamma_{ij}\ge 0$ for $i,j\neq 0$. Mean square stability of the
system is enforced by appropriate scaling of the off-diagonal entries
$\gamma_{ij}$. In these experiments, whose results are presented in Table \ref{table3},  only Algorithm 1 is considered. 

\begin{table}[h]
\caption{Application of Algorithm 1 to random examples with nonsymmetric $A_j$ and varying $\Gamma$} \label{table3}
  \begin{center}
    \begin{tabular}{|c|c|c|c|c|} \hline
         $n$&$N$   & Time (seconds) & Iteration & Residual \\ \hline 
      200 &20& 11.7 & 6.5 & 1.1e$-$07\\\hline
	100&100 & 11.3& 6 & 2.3e$-$08\\\hline
200&100 & 64.0& 7 & 1.5e$-$07\\\hline
1000&3&241&6&3.7e$-$06\\\hline
    \end{tabular}
  \end{center}
\end{table}


\subsection{A communication network example}
\label{sec:WLAN}
It has been argued in the literature (e.g.\ \cite{PlopKawk04, KawkAlle09}) that Markov jump linear systems
can be used to model communication phenomena in cyber-physical systems.
In the following we sketch such an example, where we have a fixed
number $\nu$ of entities or agents to be controlled. Each entity is
seen as a transmitting station that transmits its observed values via a
medium which follows the Carrier Sense Multiple 
Access/Collision Avoidance (CSMA/CS) principle, as it is used e.g.\ in
a WLAN transmission. A Markov jump linear system can then be used to
determine which station is allowed to send its data. 

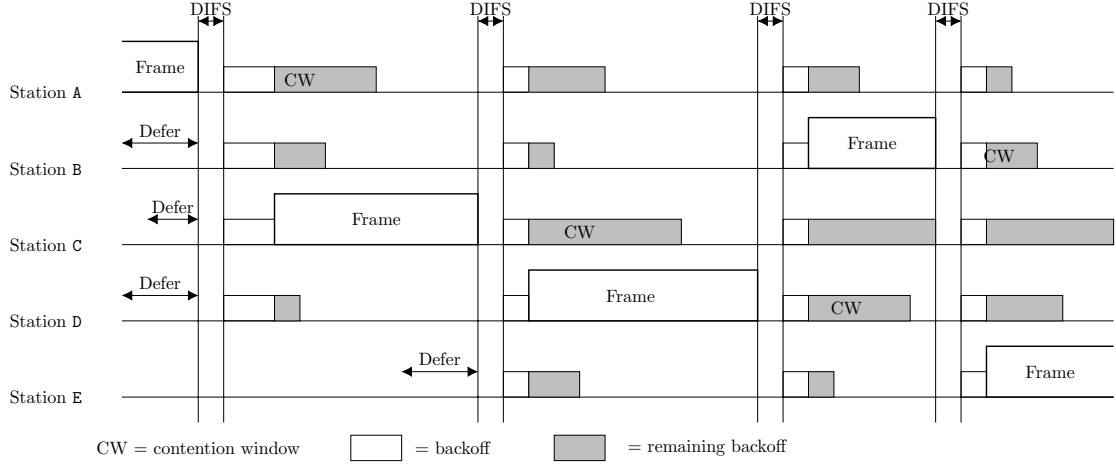
\begin{figure}[ht]
	\centering
	\resizebox{\textwidth}{!}{
	\begin{tikzpicture} \draw
(-1.5,0) 		-- (18,0)
(-1.5,-1.5) -- (18,-1.5)
(-1.5,-3) 	-- (18,-3)
(-1.5,-4.5)	-- (18,-4.5)
(-1.5,-6) 	-- (18,-6)

(0,1.5) 		-- (0,-6.5)
(0.5,1.5) 	-- (0.5,-6.5)
(5.5,1.5) 	-- (5.5,-6.5)
(6,1.5) 		-- (6,-6.5)
(11,1.5) 		-- (11,-6.5)
(11.5,1.5) 	-- (11.5,-6.5)
(14.5,1.5) 	-- (14.5,-6.5)
(15,1.5) 		-- (15,-6.5)
;

\draw[thick]													(-1.5,1) 		-- (0,1) 				-- (0,0)			 	-- (-1.5,0);

\draw 																(.5,0) 			-- (.5,.5) 			-- (1.5,.5) 		-- (1.5,0) 			-- cycle;
\filldraw[fill={rgb:black,1;white,3}]	(1.5,0) 		-- (1.5,0.5)		-- (3.5,0.5)		-- (3.5,0) 			-- cycle;
\draw 																(.5,-1.5) 	-- (.5,-1) 			-- (1.5,-1) 		-- (1.5,-1.5)		-- cycle;
\filldraw[fill={rgb:black,1;white,3}]	(1.5,-1.5) 	-- (1.5,-1) 		-- (2.5,-1) 		-- (2.5,-1.5) 	-- cycle;
\draw 																(.5,-3) 		-- (.5,-2.5) 		-- (1.5,-2.5) 	-- (1.5,-3) 		-- cycle;
\draw[thick]													(1.5,-3) 		-- (5.5,-3) 		-- (5.5,-2) 		-- (1.5,-2) 		-- cycle;
\draw 																(.5,-4.5) 	-- (.5,-4) 			-- (1.5,-4) 		-- (1.5,-4.5) 	-- cycle;
\filldraw[fill={rgb:black,1;white,3}]	(1.5,-4.5) 	-- (1.5,-4) 		-- (2,-4) 			-- (2,-4.5) 		-- cycle;

\draw 																(6,0) 			-- (6,.5) 			-- (6.5,.5) 		-- (6.5,0) 			-- cycle;
\filldraw[fill={rgb:black,1;white,3}]	(6.5,0) 		-- (6.5,0.5) 		-- (8,0.5) 			-- (8,0) 				-- cycle;
\draw 																(6,-1.5) 		-- (6,-1) 			-- (6.5,-1) 		-- (6.5,-1.5) 	-- cycle;
\filldraw[fill={rgb:black,1;white,3}]	(6.5,-1.5) 	-- (6.5,-1) 		-- (7,-1)		 		-- (7,-1.5) 		-- cycle;
\draw 																(6,-3) 			-- (6,-2.5) 		-- (6.5,-2.5) 	-- (6.5,-3) 		-- cycle;
\filldraw[fill={rgb:black,1;white,3}]	(6.5,-2.5) 	-- (6.5,-3) 		-- (9.5,-3) 		-- (9.5,-2.5) 	-- cycle;
\draw 																(6,-4.5) 		-- (6,-4) 			-- (6.5,-4) 		-- (6.5,-4.5) 	-- cycle;
\draw[thick]													(6.5,-4.5) 	-- (6.5,-3.5) 	-- (11,-3.5) 		-- (11,-4.5) 		-- cycle;
\draw 																(6,-6) 			-- (6,-5.5) 		-- (6.5,-5.5) 	-- (6.5,-6) 		-- cycle;
\filldraw[fill={rgb:black,1;white,3}]	(6.5,-5.5) 	-- (6.5,-6) 		-- (7.5,-6) 		-- (7.5,-5.5) 	-- cycle;

\draw 																(11.5,0) 		-- (11.5,.5) 		-- (12,.5) 			-- (12,0) 			-- cycle;
\filldraw[fill={rgb:black,1;white,3}]	(12,0) 			-- (12,0.5)	 		-- (13,0.5) 		-- (13,0)	 			-- cycle;
\draw 																(11.5,-1.5) -- (11.5,-1) 		-- (12,-1)			-- (12,-1.5) 		-- cycle;
\draw[thick]													(12,-1.5) 	-- (12,-0.5) 		-- (14.5,-0.5) 	-- (14.5,-1.5) 	-- cycle;
\draw 																(11.5,-3) 	-- (11.5,-2.5) 	-- (12,-2.5) 		-- (12,-3) 			-- cycle;
\filldraw[fill={rgb:black,1;white,3}]	(12,-2.5) 	-- (12,-3) 			-- (14.5,-3) 		-- (14.5,-2.5) 	-- cycle;
\draw 																(11.5,-4.5) -- (11.5,-4) 		-- (12,-4) 			-- (12,-4.5) 		-- cycle;
\filldraw[fill={rgb:black,1;white,3}]	(12,-4.5) 	-- (12,-4) 			-- (14,-4) 			-- (14,-4.5) 		-- cycle;
\draw 																(11.5,-6) 	-- (11.5,-5.5) 	-- (12,-5.5) 		-- (12,-6) 			-- cycle;
\filldraw[fill={rgb:black,1;white,3}]	(12,-6) 		-- (12,-5.5) 		-- (12.5,-5.5) 	-- (12.5,-6) 		-- cycle;

\draw 																(15,0) 			-- (15,.5) 			-- (15.5,.5) 		-- (15.5,0) 		-- cycle;
\filldraw[fill={rgb:black,1;white,3}]	(15.5,0) 		-- (15.5,0.5) 	-- (16,0.5) 		-- (16,0) 			-- cycle;
\draw 																(15,-1.5) 	-- (15,-1) 			-- (15.5,-1) 		-- (15.5,-1.5) 	-- cycle;
\filldraw[fill={rgb:black,1;white,3}]	(15.5,-1.5) -- (15.5,-1) 		-- (16.5,-1) 		-- (16.5,-1.5) 	-- cycle;
\draw 																(15,-3) 		-- (15,-2.5) 		-- (15.5,-2.5) 	-- (15.5,-3) 		-- cycle;
\filldraw[fill={rgb:black,1;white,3}]	(15.5,-2.5) -- (15.5,-3) 		-- (18,-3) 			-- (18,-2.5) 		-- cycle;
\draw 																(15,-4.5) 	-- (15,-4) 			-- (15.5,-4) 		-- (15.5,-4.5) 	-- cycle;
\filldraw[fill={rgb:black,1;white,3}]	(15.5,-4.5) -- (15.5,-4) 		-- (17,-4) 			-- (17,-4.5) 		-- cycle;
\draw 																(15,-6) 		-- (15,-5.5) 		-- (15.5,-5.5) 	-- (15.5,-6) 		-- cycle;
\draw[thick]													(18,-5)			-- (15.5,-5) 		-- (15.5,-6) 		-- (18,-6);

\draw 																(3,-6.75) 	-- (3,-7.25) 		-- (4,-7.25) 		-- (4,-6.75) 		-- cycle;
\filldraw[fill={rgb:black,1;white,3}]	(7,-6.75) 	-- (7,-7.25) 		-- (8,-7.25) 		-- (8,-6.75) 		-- cycle;

\draw[<->] (0,1.4) 		-- node[midway,above] {DIFS} 	(0.5,1.4);
\draw[<->] (5.5,1.4) 	-- node[midway,above] {DIFS} 	(6,1.4);
\draw[<->] (11,1.4) 	-- node[midway,above] {DIFS} 	(11.5,1.4);
\draw[<->] (14.5,1.4) -- node[midway,above] {DIFS} 	(15,1.4);

\draw[<->] (-1.5,-1) 	-- node[midway,above] {Defer} (0,-1);
\draw[<->] (-1,-2.5) 	-- node[midway,above] {Defer} (0,-2.5);
\draw[<->] (-1.5,-4) 	-- node[midway,above] {Defer} (0,-4);
\draw[<->] (4,-5.5) 	-- node[midway,above] {Defer} (5.5,-5.5);

\draw
(-3,0)				node[] (A) 			{Station \texttt{A}}
(-3,-1.5)			node[] (B) 			{Station \texttt{B}}
(-3,-3)				node[] (C) 			{Station \texttt{C}}
(-3,-4.5)			node[] (D) 			{Station \texttt{D}}
(-3,-6)				node[] (E) 			{Station \texttt{E}}

(-0.75,0.5)		node[] (Frame1)	{Frame}
(3.5,-2.5)		node[] (Frame2) {Frame}
(8.5,-4)			node[] (Frame3) {Frame}
(13.25,-1)		node[] (Frame4) {Frame}
(16.75,-5.5)	node[] (Frame5) {Frame}

(2,.25) 			node[] (CW1)		{CW}
(7.5,-2.75) 	node[] (CW2)		{CW}
(12.75,-4.25)	node[] (CW3)		{CW}
(15.75,-1.25)	node[] (CW3)		{CW}

(0,-7)				node[]	(CW)		{CW = contention window}
(5,-7)				node[]	(bkoff)	{= backoff}
(10,-7)				node[]	(remb)	{= remaining backoff}
;
\end{tikzpicture}
}
	\caption{Backoff procedure: We can see 5 stations, where
          station \texttt{A} sends in the beginning. During the transmission of
          station \texttt{A} the access  of \texttt{B}, \texttt{C} and \texttt{D} to the physical medium is
          deferred until the medium is idle for a certain amount of
          time, called a \emph{{\bf Di}stributed Coordination {\bf F}unction
            Interframe {\bf S}pace Period}  (DIFS). Station \texttt{A} then gets a new
          backoff value in the according contention window. While the
          physical medium is idle, the backoff value of all
          ready-to-transmit stations
         is reduced until one of them reaches zero, here station \texttt{C}. Station \texttt{C} then begins transmitting and the procedure repeats itself.
The probability for a station to be the next to transmit
depends on the remaining backoff value and thus on the time since its
last transmission. 
}
	\label{Backoff}
\end{figure}

In a CSMA/CA network after every transmission a backoff value is
assigned to the station which did just transmit. The backoff value is a
uniformly distributed random variable in a given interval
\texttt{[0,ContentionWindowMax]}. While the physical medium of a
CSMA/CS network is idle, the backoff value of each station is reduced
continuously. When the backoff value reaches $0$ the 
respective station is again allowed to transmit its data. See figure
\ref{Backoff} for a visualization.

We then keep track of the $\tau$ last transmissions 
and whether the last
transmission was faulty or not. 
This information can then be encoded as states in a transition rate
matrix for a Markov chain. Assuming $2$ entities and a transition rate
matrix which keeps track of the last $3$ transmissions this leads to
a transition rate matrix of size $16 \times 16$, which we present here in Table \ref{table:Pi}.

\begin{table}[h]
\resizebox{\textwidth}{!}{\begin{minipage}{\textwidth}
  
\begin{tabular}{@{\thinspace}l @{\thinspace\thinspace\thinspace}l @{\thinspace}c @{\thinspace}l@{\thinspace} | @{\thinspace}c @{\thinspace}c @{\thinspace}c @{\thinspace}c @{\thinspace}c @{\thinspace}c @{\thinspace}c @{\thinspace}c @{\thinspace}c @{\thinspace}c @{\thinspace}c @{\thinspace}c @{\thinspace}c @{\thinspace}c @{\thinspace}c @{\thinspace}c}
&&	& e& \texttt{0}& \texttt{0}& \texttt{0}& \texttt{0}& \texttt{0}& \texttt{0}& \texttt{0}& \texttt{0}& \texttt{1}& \texttt{1}& \texttt{1}& \texttt{1}& \texttt{1}& \texttt{1}& \texttt{1}& \texttt{1} \\
&&	& t+1& \texttt{A}& \texttt{A}& \texttt{A}& \texttt{A}& \texttt{B}& \texttt{B}& \texttt{B}& \texttt{B}& \texttt{A}& \texttt{A}& \texttt{A}& \texttt{A}& \texttt{B}& \texttt{B}& \texttt{B}& \texttt{B} \\
&&	& t& \texttt{A}& \texttt{A}& \texttt{B}& \texttt{B}& \texttt{A}& \texttt{A}& \texttt{B}& \texttt{B}& \texttt{A}& \texttt{A}& \texttt{B}& \texttt{B}& \texttt{A}& \texttt{A}& \texttt{B}& \texttt{B} \\ 
&&	& t-1& \texttt{A}& \texttt{B}& \texttt{A}& \texttt{B}& \texttt{A}& \texttt{B}& \texttt{A}& \texttt{B}& \texttt{A}& \texttt{B}& \texttt{A}& \texttt{B}& \texttt{A}& \texttt{B}& \texttt{A}& \texttt{B} \\
\hline
e& t& t-1& t-2&  &  &  &  &  &  &  &  &  &  &  &  &  &  &  &\\
\hline
\texttt{0}& \texttt{A}& \texttt{A}& \texttt{A}& 0.19& 0  & 0  &0 &0.78& 0  & 0  & 0  & 0.01& 0  & 0  & 0  & 0.02& 0  & 0  & 0\\
\texttt{0}& \texttt{A}& \texttt{A}& \texttt{B}& 0.32& 0  & 0  & 0  & 0.65& 0  & 0  & 0  & 0.01  &0  & 0  & 0  & 0.02& 0  & 0  & 0\\
\texttt{0}& \texttt{A}& \texttt{B}& \texttt{A}& 0  &0.42& 0  & 0  & 0  & 0.55& 0  & 0  & 0  & 0.01& 0  & 0  & 0  & 0.02& 0  & 0\\ 
\texttt{0}& \texttt{A}& \texttt{B}& \texttt{B}& 0  &0.42& 0  & 0  & 0  & 0.55& 0  & 0  & 0  &0.01& 0  & 0  & 0  & 0.02& 0  & 0\\
\texttt{0}& \texttt{B}& \texttt{A}& \texttt{A}& 0  & 0  &0.55& 0  & 0  & 0  & 0.42& 0  & 0  & 0  & 0.02& 0  & 0  & 0  & 0.01& 0\\
\texttt{0}& \texttt{B}& \texttt{A}& \texttt{B}& 0  & 0  &0.55& 0  & 0  & 0  & 0.42  & 0& 0  & 0  & 0.02  & 0& 0  & 0  & 0.01  & 0\\
\texttt{0}& \texttt{B}& \texttt{B}& \texttt{A}& 0  & 0  & 0  &0.65& 0  & 0  & 0  & 0.32& 0  & 0  & 0  & 0.02& 0  & 0  & 0  & 0.01\\
\texttt{0}& \texttt{B}& \texttt{B}& \texttt{B}& 0  & 0  & 0  &0.78& 0  & 0  & 0  & 0.19& 0  & 0  & 0  & 0.02& 0  & 0  & 0  & 0.01\\
\texttt{1}& \texttt{A}& \texttt{A}& \texttt{A}& 0.25& 0  & 0  & 0  & 0  &0  & 0  & 0  & 0.75& 0  & 0  & 0  & 0  &0  & 0  & 0  \\
\texttt{1}& \texttt{A}& \texttt{A}& \texttt{B}& 0.25& 0  & 0  & 0  & 0  &0  & 0  & 0  & 0.75& 0  & 0  & 0  & 0  &0  & 0  & 0  \\
\texttt{1}& \texttt{A}& \texttt{B}& \texttt{A}& 0  &0.25& 0  & 0  & 0  & 0  & 0  & 0  & 0  &0.75& 0  & 0  & 0  & 0  & 0  & 0  \\
\texttt{1}& \texttt{A}& \texttt{B}& \texttt{B}& 0  &0.25& 0  & 0  & 0  & 0  & 0  & 0  & 0  &0.75& 0  & 0  & 0  & 0  & 0  & 0  \\
\texttt{1}& \texttt{B}& \texttt{A}& \texttt{A}& 0  & 0  & 0  & 0  & 0  & 0  &0.25& 0  & 0  & 0  & 0  & 0  & 0  & 0  &0.75& 0  \\
\texttt{1}& \texttt{B}& \texttt{A}& \texttt{B}& 0  & 0  & 0  & 0  & 0  & 0  &0.25& 0  & 0  & 0  & 0  & 0  & 0  & 0  &0.75& 0  \\
\texttt{1}& \texttt{B}& \texttt{B}& \texttt{A}& 0  & 0  & 0  & 0  & 0  & 0  & 0  & 0.25& 0  & 0  & 0  & 0  & 0  & 0  & 0  & 0.75\\
\texttt{1}& \texttt{B}& \texttt{B}& \texttt{B}& 0  & 0  & 0  & 0  & 0  & 0  & 0  & 0.25& 0  & 0  & 0  & 0  & 0  & 0  & 0  & 0.75	\\
\end{tabular}
\caption[]{Transition matrix $\Theta$ for $\tau=3$, $\nu=2$ with
  stations \texttt{A} and \texttt{B}.}
\label{table:Pi}
\end{minipage} }
\end{table}
In the matrix we see the probabilities to get from the current
situation, presented in the row of the matrix, to the future situation
presented in the columns. We denote by $e = 0$ that the transmission
is correct and $e=1$ that the transmission is incorrect. The values
for $t$ denote which station sends at the moment, $t-1$ which station
did send before and $t-2$ which station did send before $t-1$. The
station sending next is then denoted by $t+1$. The
probability for an error in the next transmission is set to $0.03$, if
the current transmission is correct, and $0.75$ if we have an
error in the current transmission, as e.g.\ suggested in \cite{Gilb60}, see also \cite
{Burst}.
 
Now we want to approximate the probability for a fixed station to be
 the next to send. 
Let
$$\bar{w} = \Big(\nu -
  \left\vert J\right\vert + \sum_{j\in J}
  \frac{1}{\tau+1-\left(j-1\right)}\Big)^{-1}\;,$$ 
where for each station, the set $J$ contains the index of its last
occurrence in the memory.
Then the approximated probability for a fixed station to be the next to send
(neglecting the probability of an error), is either $\bar w$, if it
did not send in the last $\tau$ transmissions, or
$\frac{\bar{w}}{\tau+1-\left(i-1\right)}$, if  its last transmission
was at $t+1-i$. The transition matrix $\Theta$ in Table \ref{table:Pi} contains
these entries. As an example consider the transition \texttt{0AAA}
$\leadsto$ \texttt{0BAA}.
Here $J=\{1\}$ because the first position of \texttt{AAA} contains the station
\texttt{A}. Hence $\bar w=(2-1+\frac1{3+1-(1-1)})^{-1}=\frac45$. Thus the
transition probability is  $0.8\cdot0.97\approx 0.78=\Theta_{1,5}$.
For a different case consider \texttt{0BBA} $\leadsto$ \texttt{0ABB}. Here $J=\{1,3\}$
because the first occurences of \texttt{B} and \texttt{A} in \texttt{BBA} are at positions $1$
and $3$. Here $\bar
w=(2-2+\frac1{3+1-(1-1)}+\frac1{3+1-(3-1)})^{-1}=\frac43$ and the
transition probability is $\frac43\cdot
\frac12\cdot0.97\approx0.65=\Theta_{7,4}$. 
The formula was obtained by a discretization of the interval
 \texttt{[0,ContentionWindow]} and an approximation for the exact
 probability for each station to send next.  Details of the technical derivation
 can be found in \cite{Vier17}. In the case of an error, the
 transmission memory is not changed until the error state changes. So
 for arbitrary \texttt{X,Y,Z}, only \texttt{1XYZ}$\leadsto$\texttt{1XYZ} or \texttt{1XYZ}$\leadsto$\texttt{0XYZ}  is possible, with the corresponding probabilities.
Note that we have thus approximated  the transition probabilities for a discrete time
setting. 
The transformation to the continuous time
situation is obtained by  
$\Gamma= a(\Theta - I)$ where $a$ is the average time spent in one
transmission mode, see
\cite{CTMC}. 

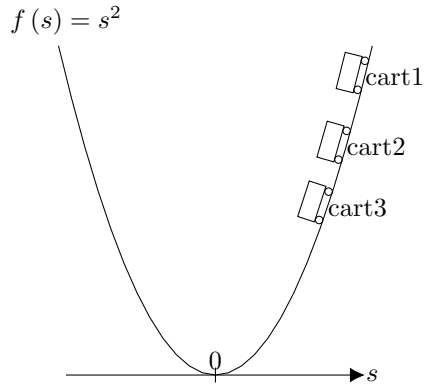
\begin{figure}[ht]
	\centering
	\resizebox{.4\textwidth}{!}{
	\begin{tikzpicture}
\draw [domain=-2.1:2.1] plot (\x, {\x*\x});

\draw
(2+0.0485-0.0485,4.1964+0.0121) node[draw,circle,inner sep=0pt,minimum size=.1cm] (wheel11){}
(2-0.0485-0.0485,3.8083+0.0121) node[draw,circle,inner sep=0pt,minimum size=.1cm] (wheel12){}
(2-0.21,4+0.051) node[draw, rectangle, minimum width=0.2cm,minimum height=0.5cm, rotate=-14.04] (cart1) {}

(2+0.4,4) node[] (text1){ cart1}

(1.75+0.0549-0.0485,3.2578+0.0121) node[draw,circle,inner sep=0pt,minimum size=.1cm] (wheel21){}
(1.75-0.0549-0.0485,2.8732+0.0121) node[draw,circle,inner sep=0pt,minimum size=.1cm] (wheel22){}
(1.75-0.21,3.0625+0.065) node[draw, rectangle, minimum width=0.2cm,minimum height=0.5cm, rotate=-15.95] (cart2) {}

(1.75+0.4,3.0625) node[] (text2){cart2}

(1.5+0.0632-0.0485,2.4437+0.0121) node[draw,circle,inner sep=0pt,minimum size=.1cm] (wheel31){}
(1.5-0.0632-0.0485,2.0643+0.0121) node[draw,circle,inner sep=0pt,minimum size=.1cm] (wheel32){}
(1.5-0.21,2.25+0.075) node[draw, rectangle, minimum width=0.2cm,minimum height=0.5cm, rotate=-18.0] (cart3) {}

(1.5+0.4,2.25) node[] (text3){ cart3}

(-2,4.75) node[] (text1){$f\left(s\right) = s^{2}$}
(2.1,0) node[] (text1){$s$}
(0,0.2) node[] (text1){ $0$}
;

\draw [->] (-2,0) -- (2,0);
\draw (0,-0.1) -- (0,0.1);

\end{tikzpicture}
}
	\caption{A system of three carts on a parabula shaped track.
          Each cart is steered individually and can transmit its
          position $s(t)$ and velocity $v(t)$. The origin
          $(s,v)=(0,0)$ is an asymptotically stable equilibrium for all carts.}
	\label{system}
\end{figure}

Together with this model  of communication, we consider the simple system described in figure \ref{system}.

We denote with $e\in \left\{0,1\right\}$ whether or not an error
ocurred in the last transmission and with $j\in\left\{1,2,3\right\}$
the station which did send last and its corresponding diagonal block.
Note that both $e$ and $j$ are determined by $r\left(t\right)$ but are
added here for easier understanding of the model.
The state vector of the whole system is $\begin{pmatrix}
		{s_{1}\left(t\right)} & {v_{1}\left(t\right)} & {s_{2}\left(t\right)} & {v_{2}\left(t\right)} & {s_{3}\left(t\right)} & {v_{3}\left(t\right)}
		\end{pmatrix}^{T}$. 
The linearized dynamics  are given by the equations
\begin{align*} \label{car_model}
		\dot{x}\left(t\right) &= Ax\left(t\right) + Bu\left(t\right)\;,\quad 
				y\left(t\right) = C\left(r\left(t\right)\right)x\left(t\right)\;,\quad\text{ where}\\		
				A &= \begin{bmatrix}
		\bar{A} & & \\
		& \bar{A} & \\
		& & \bar{A} \\		
 \end{bmatrix}\;,\quad
				\quad \bar{A} = \begin{bmatrix}
0 & 1 \\
-mgs
& -R \\
			 \end{bmatrix}\;,\quad
				B = \begin{bmatrix}
		\bar{B} & & \\
		& \bar{B} & \\
		& & \bar{B} \\		
 \end{bmatrix}\;,\quad
				\quad \bar{B} = \begin{bmatrix}
	0 \\
	1 \\
			 \end{bmatrix}\;, \\
				C_{1} &= \begin{bmatrix}
	(1-e)I_{2} & 0_{2} & 0_{2}	
\end{bmatrix}, \quad C_{2} = \begin{bmatrix}
	0_{2} & (1-e)I_{2} & 0_{2}
\end{bmatrix}, \quad C_{3} = \begin{bmatrix}
	0_{2} & 0_{2} & (1-e)I_{2} 
\end{bmatrix}\;.
\end{align*}
Here $I_{2}$ denotes the identity matrix and  $0_{2}$ the $0$ matrix in
$\mathbb{R}^{2\times 2}$, while $m,g,R$ denote the mass of each cart,
the gravitational acceleration, and the friction coefficient,
respectively.
In this application only the matrices $C_{i}$ are depending on the
Markov process. There is still a need for fast solving of Lyapunov
like equations in this applications since the size of the transition
rate matrix, and therefore the number of coupled Lyapunov equations,
scales with a factor $\nu^{\tau}$, where $\nu$ is the number of entities
and $\tau$ is the number of tracked transmissions. 

We compute the observability Gramian \eqref{lya2}.
Let $\nu=3$, $\tau=3$, $m=1$, $g=9.81$, and $R=0.1$.
Table~\ref{table6} shows the computational time, iteration number, and
the residual in Algorithm 1. The results show that the algorithm works well in
practical applications with relevant dimensions.

\begin{table}[h]
\caption{Application of Algorithm 1 to networked control system} \label{table6}
  \begin{center}
    \begin{tabular}{|c|c|c|c|c|c|} \hline
         $\nu$& $n$&$N$  & Time (seconds) & Iteration & Residual \\ \hline 
    3&  6&54& 0.04 & 3.5 & 2.3e$-$11 \\ \hline
5&10&250&1.1 &13.5&2.2e$-$12    \\ \hline
10&20&2000&27.3  &14  &1.5e$-$10  \\ \hline
    \end{tabular}
  \end{center}
\end{table}


\section{Conclusion} \label{sec6}

We have compared optimization based methods and a preconditioned Krylov
subspace iteration for the solution of Lyapunov equations related to
Markov jump linear systems. From our numerical experiments we deduce
that only the Krylov subspace iteration lends itself for large
systems. As an application we have sketched a Markov jump linear system
model for a networked control system with WLAN based communication. In
ongoing research, we plan to elaborate further on this model. We
expect that efficient Lyapunov solvers will be an essential tool. A
next step in this direction could be the development of low-rank methods
like they were considered e.g.\ in \cite{BennBrei13, ShanSimo16} for
other types of Lyapunov equations.




\bibliographystyle{abbrv}
\bibliography{NACO.bib}

\medskip

\medskip
\end{document}